\documentclass[11pt]{article}
\usepackage[T2A]{fontenc}
\usepackage[cp1251]{inputenc}

\usepackage[T2A]{fontenc}
\usepackage[cp1251]{inputenc}
\usepackage{amsfonts}
\usepackage{enumerate}
\usepackage{amssymb}
\usepackage{amsmath}

\begin{document}

\begin{center}{\Large \textbf{On a characterization of idempotent
distributions on discrete fields and on the field of $p$-adic
numbers}}
\end{center}

\bigskip

\centerline {\bf G.M. Feldman and M.V. Myronyuk}


\bigskip

\begin{abstract}
We prove the following theorem. Let $X$ be a discrete field,  $\xi$
and $\eta$ be independent identically distributed random variables
with values in $X$ and distribution  $\mu$. The random variables
$S=\xi+\eta$ and
 $D=(\xi-\eta)^2$ are independent if and only if
 $\mu$ is an idempotent distribution.
A similar result is also proved in the case when $\xi$ and $\eta$
are independent identically distributed random variables with values
in the field of  $p$-adic numbers  $\mathbf{Q}_p$, where $p>2$,
assuming that the distribution $\mu$ has a continuous density.
\end{abstract}

{\it  Keys words and phrases:} Characterization theorem, idempotent
distribution, discrete field, the field of  $p$-adic numbers.


\section{Introduction}
\label{intro}

According to the classical Kac--Bernstein theorem if $\xi$ and
$\eta$ are independent random variables and their sum  $\xi +\eta$
and difference
  $\xi -\eta$ are independent, then the random variables $\xi$ and
$\eta$ are Gaussian. This characterization of a Gaussian measure
remains true if instead of the difference $\xi -\eta$ we consider
its square
 $(\xi -\eta)^2,$ but assume that the random variables
  $\xi$ are  $\eta$ identically distributed. This characterization of a Gaussian
  measure is a particular case of the well-known Geary theorem: if  $\xi_1, \xi_2,
\dots, \xi_n, \ n\ge 2$, are independent identically distributed
random variables such that the sample mean  ${\bar \xi}={1\over
n}\sum_{j=1}^n\xi_j$ and
 the sample variance $s^2={1\over n}\sum_{j=1}^n(\xi_j-{\bar \xi})^2$
are independent, then the random variables $\xi_j$ are Gaussian
(\cite{Ge}, \cite{Lu}, \cite{KaSa}, \cite{Zi}).

A lot of research are devoted to generalizations of the classical
characterization theorems of mathematical statistics to different
algebraic structures,
 first of all to locally compact Abelian groups (see e.g.
 \cite{Fe10}--\cite{Fe19}, \cite{Fe30}, \cite{JOTP}, \cite{GraLo},  \cite{M}--\cite{NeueScho},   and also \cite{Fe1},
 where one can find additional references). So far they dealt only with linear forms of independent random variables with
values in a group, where coefficients of the linear forms
 were topological automorphisms of the group (the Kac--Bernstein theorem, the Skitovich--Darmois theorem, etc).

In this article we study the simplest "nonlinear" characterization
problem. Let  $X$ be a locally compact field. Let  $\xi$ and $\eta$
be independent identically distributed random variables with values
in $X$ and  distribution  $\mu$. Assume that the random variables
$S=\xi +\eta$ and $D=(\xi -\eta)^2$ are independent. Assuming that
$X$ is a discrete field, we prove that then $\mu$ is idempotent
distribution. A similar result we prove in the case, when $X$ is the
field of $p$-adic numbers $\mathbf{Q}_p$, where $p>2$, assuming that
$\mu$ has a density with respect to a Haar measure on
$\mathbf{Q}_p$, and this density is continuous. Since on totally
disconnected locally compact Abelian groups, in particular, on
discrete groups and on $\mathbf{Q}_p$ idempotent distributions are a
natural analogue of Gaussian measures, the obtained results one can
consider as  a natural analogue of the corresponding
characterization theorem both for discrete fields and the field
$\mathbf{Q}_p$.

Let  $X$ be a locally compact field.
 Denote by $e$ the unit of the field $X$. The additive group of the field  $X$
 is a locally compact Abelian group. We also denote this group by $X$.
   Denote by $(x,y),$ $x, y\in X,$
elements of the group $X^2$. Denote by $T$ the mapping $T:X^2\mapsto
X^2$ defined by the formula $T(x, y)=(x+y, (x-y)^2)$. Let $A$ be a
subset of  $X$. Put $A^{[2]}=\{x\in X:x=t^2, \ t\in A\}$.  Let $G$
be an arbitrary locally compact Abelian group. Denote by $m_G$ a
Haar measure on $G$. If $G$ is a compact group, then we assume that
$m_G$ is a distribution, i.e.
 $m_G(G)=1$. Denote by $I(X)$ the set of all idempotent distributions on
 $X$, i.e. the set of all shifts of Haar distributions $m_K$ of compact
  subgroup $K$ of the group $X$.
If $x\in X$ denote by $E_x$ the degenerate distribution concentrated
at the point
 $x$. If $\xi$ and $\eta$ are random variables with values in
  $X$, then we denote by $\mu_{\xi}$  the distribution of the random
  variable  $\xi$, and by
$\mu_{(\xi, \eta)}$ the distribution of the random vector $ (\xi,
\eta)$.

\section{Characterization of idempotent distributions on discrete fields}
\label{sec:1}

We study in this section the case when independent random variables
take values in a discrete field $X$. First we consider the case when
the characteristic of the field $X$ differs from 0, and then the
case  when the characteristic of the field $X$ is equal to 0.

Note that if the characteristic of a field $X$ differs from 2, then
$X$ is an Abelian group  with unique division by  2, i.e. for any
element $x\in X$
 the element ${x\over
2}\in X$ is uniquely determined.

For proving theorems characterizing  idempotent distributions on
discrete fields we need two lemmas.

\textbf{Lemma 1.} \textit{Let $X$ be a countable discrete field of
characteristic $p$, where $p\ne 2$. Let $\xi$ and $\eta$ be
independent identically distributed random variables with values in
$X$ and distribution
 $\mu$  such that $\mu(0)>0$. In order that the random variables
  $S=\xi+\eta$ and $D=(\xi-\eta)^2$ be independent it is necessary and sufficient that
  the function $\mu(x)$  satisfy the equation}
\begin{equation}\label{l1.1.}
    \mu^2(u)\mu(v)\mu(-v)=\mu^2(0) \mu(u+v)\mu(u-v), \quad u,v\in
    X.
\end{equation}

\textbf{Proof.}
 Note that $\mu_{(S,D)}=T(\mu_{(\xi,\eta)})$.
It is obvious that
\begin{equation}\label{l1.10}
    \mu_{(S,D)}(u,0)=T(\mu_{(\xi,\eta)})(u,0)=
    \mu_{(\xi,\eta)}(T^{-1}(u,0))=\mu_{(\xi,\eta)}\left({u\over
    2},{u\over 2}\right).
\end{equation}
Let    $v={t}^2, \ {t}\in X, \ {t}\ne 0$.  We have
\begin{equation}\label{l1.2}
    \mu_{(S,D)}(u,v)=\mu_{(S,D)}(u,{t}^2)=
    \mu_{(\xi,\eta)}(T^{-1}(u,{t}^2))=$$ $$
    =\mu_{(\xi,\eta)}\left({u+{t}\over 2}, {u-{t}\over 2}\right)+
    \mu_{(\xi,\eta)}\left({u-{t}\over 2},{u+{t}\over
    2}\right), \quad u\in X, \ {t}\neq 0.
\end{equation}
Taking into account the independence of the random variables  $\xi$
and $\eta$, it follows from (\ref{l1.10}) and (\ref{l1.2}) that
\begin{equation}\label{l1.11}
     \mu_{(S,D)}(u,0) =  \mu^2\left({u\over 2}\right), \quad u\in X,
\end{equation}
\begin{equation}\label{l1.12}
    \mu_{(S,D)}(u, v) = 2\mu\left({u+t\over 2}\right)\mu\left({u-t\over 2}\right),
    \quad u\in X, \ v=t^2, \ t\neq 0.
\end{equation}

 First we prove the necessity. Taking into account the independence of the random
 variables $S$ and $D$, it follows from
(\ref{l1.11}) and (\ref{l1.12}) that
\begin{equation}\label{l1.11ab}
     \mu_S(u) \mu_D(0) =  \mu^2\left({u\over 2}\right), \quad u\in X,
\end{equation}
\begin{equation}\label{l1.12ab}
   \mu_S(u) \mu_D(t^2) = 2\mu\left({u+t\over 2}\right)\mu\left({u-t\over 2}\right),
   \quad u\in X, \ t\neq 0.
\end{equation}
Moreover (\ref{l1.11ab}) implies that $\mu_D(0)\neq 0$, because
otherwise $\mu(x)=0$ for all $x\in X$, contrary to the condition of
the theorem. Find $\mu_S(u)$ from (\ref{l1.11ab}) and substitute it
to
 (\ref{l1.12ab}). We get
\begin{equation}\label{l1.13}
     \mu^2\left({u\over 2}\right)\mu_D(t^2)=
     2\mu_D(0)\mu\left({u+t\over 2}\right)\mu\left({u-t\over 2}\right),
     \quad u\in X, \ t\neq 0.
\end{equation}
Rewrite   (\ref{l1.13}) in the form
\begin{equation}\label{l1.131}
     \mu^2(u)\mu_D(4t^2)=2\mu_D(0)\left(\mu\left({u+t}\right)
     \mu\left({u-t}\right)\right),
     \quad u\in X, \ t\neq 0.
\end{equation}
Putting in (\ref{l1.131}) $u=0$, we obtain
\begin{equation}\label{l1.14}
     \mu_D(4t^2)={2\mu_D(0)\mu(t)\mu(-t)\over \mu^2(0)},
     \quad t\neq 0.
\end{equation}
It follows from (\ref{l1.131}) and that (\ref{l1.14})
\begin{equation}\label{l1.15}
     \mu^2(u)\mu(t)\mu(-t)=\mu^2(0)\mu(u+t)\mu(u-t),
     \  u\in X, \ t\neq 0.
\end{equation}
Note that (\ref{l1.15}), obviously, holds true also when $t=0$.
Therefore the function $\mu(x)$   satisfies equation (\ref{l1.1.}).
So, we proved the necessity.

To prove sufficiency we first note that   random variables $\xi_1$
and $\xi_2$ with values in a discrete space  $X$ are independent if
and only if there exist functions  $a(x)$   and $b(x)$ on $X$ such
that
\begin{equation} \label{e16.01.2}
\mu_{(\xi_1, \xi_2)}(x,y)=a(x)b(y), \quad x, y\in X.
\end{equation}

Let $\xi$ and $\eta$ be independent identically distributed random
variables with values in $X$ and distribution
 $\mu$  such that $\mu(0)>0$.
It  is obvious that if $v\notin X^{[2]}$, then $\mu_{(S,D)}(u,v)=0$.
It follows from (\ref{l1.11})   and (\ref{l1.12}) that we have a
representation
$$
\mu_{(S, D)}(u, v)=
\left\{%
\begin{array}{ll}
    \mu^2\left({u\over 2}\right), & \hbox{$u\in X, v=0$;} \\
    2\mu\left({u+t\over 2}\right)\mu\left({u-t\over 2}\right), & \hbox{$u\in X, v=t^2, t\ne 0$;} \\
    0, & \hbox{$u\in X, v\notin X^{[2]}$.} \\
\end{array}%
\right.
$$
This representation implies that if a function  $\mu(x)$ satisfies
equation (\ref{l1.1.}), then the random variables   $S$ and $D$ are
independent, because in this case we can put in (\ref{e16.01.2})
$a(x)=\mu^2\left({x\over 2}\right)$ and
$$ b(x)=
\left\{%
\begin{array}{ll}
    1, & \hbox{$x=0$;} \\
    \displaystyle{{2\mu\left({t\over 2}\right)\mu\left(-{t\over 2}\right)\over
\mu^2(0)}}, & \hbox{$x=t^2,  t\ne 0$;} \\
    0, & \hbox{$x\notin X^{[2]}$.} \\
\end{array}%
\right.
$$

\textbf{Lemma 2.} \textit{Let $X$ be an Abelian group  with unique
division by  2. Let $\mu$ be a function on $X$, satisfying equation
$(\ref{l1.1.})$ and the condition $\mu(0)>0$. Then the set $K=\{x\in
X: \mu(x)\ne 0\}$ is a subgroup of $X$.}

\textbf{Proof.} Assume that $\mu(x)\ne 0$ at a point $x \in X.$ Put
in (\ref{l1.1.}) $u=v={x\over 2}$. We get
\begin{equation}\label{t1.1}
    \mu^3\left({x\over 2}\right)\mu\left(-{x\over 2}\right)=\mu^3(0) \mu(x).
\end{equation}
Since $\mu(0) \ne 0$ and $\mu(x)\ne 0$, it follows from (\ref{t1.1})
that
\begin{equation}\label{t1.2}
    \mu\left({x\over 2}\right)\mu\left(-{x\over 2}\right)\ne 0.
\end{equation}
Put in (\ref{l1.1.}) $u=v=-{x\over 2}$. We obtain
\begin{equation}\label{t1.17}
     \mu^3\left(-{x\over 2}\right)\mu\left({x\over 2}\right)=\mu^3(0) \mu(-x).
\end{equation}
Taking into account (\ref{t1.2}), we imply from (\ref{t1.17}) that
$\mu(-x)\ne 0$. So, we proved that if  $\mu(x)\ne 0$, then also
$\mu(-x)\ne 0$. When it is considered, (\ref{l1.1.}) implies that
the set $K=\{x\in X: \mu(x)\ne 0\}$ is a subgroup of $X$.

First we prove a theorem on characterization of idempotent
distributions for a discrete fields  of nonzero characteristic.

\textbf{Theorem 1.} \textit{Let $X$ be a countable discrete field of
characteristic $p$, where $p>2$. Let $\xi$ and $\eta$ be independent
identically distributed random variables with values in $X$  and
distribution
 $\mu$.  In order that the random variables
  $S=\xi+\eta$ and $D=(\xi-\eta)^2$ be independent it is necessary
  and sufficient that  $\mu\in I(X)$.}

\textbf{Proof.}  Necessity. Replacing if it is necessary the random
variables $\xi$ and $\eta$ by new independent random variables
$\xi+x$ and $\eta+x$, we can assume from the beginning that
$\mu(0)>0$. Then Lemma 1 implies that the function  $\mu(x)$
satisfies equation (\ref{l1.1.}). Since the characteristic of the
field  $X$ is greater than 2,   $X$ is an Abelian group  with unique
division by  2. Hence, by Lemma 2 the set $K=\{x\in X: \mu(x)>0\}$
is a subgroup of $X$. Let $x_0\in K$, $x_0\ne 0$. Denote by $L$ the
subgroup generating by  $x_0$. Then $L\cong\mathbf{Z}(p)$, where
$\mathbf{Z}(p)$ is the group of residue modulo $p$, and $L\subset
K$. Consider the restriction of equation (\ref{l1.1.}) to $L$. Put
$\varphi(x)=\log \mu(x)$, $x\in L$. It follows from (\ref{l1.1.})
that
$$
2 \varphi(u)+\varphi(v)+\varphi(-v)=
2\varphi(0)+\varphi(u+v)+\varphi(u-v), \quad u ,v \in L.
$$
Integrate both sides of this equality by the measure
 $dm_L(v)$. We get that
$\varphi(u)=\varphi(0)$ for all $u\in L$. It means that
$\mu(x)=\mu(0)$ for all $x\in L$, and hence, $\mu(x)=\mu(0)$ for all
$x\in K$. This implies that  $K$ is a finite subgroup and $\mu=m_K$.
The necessity is proved.

Sufficiency.  Let $K$ be a finite subgroup of $X$. Denote by $|K|$
the number of elements of the subgroup $K$. Let $\xi$ and $\eta$ be
independent identically distributed random variables with values in
$X$  and distribution  $\mu=m_K$. Since all nonzero elements of the
field $X$ have order $p$, and $p>2$, the following statement holds:
if $2x\in K,$ then $x\in K$. This easily implies that the function
$$ \mu(x)=m_K(x)=
\left\{%
\begin{array}{ll}
    |K|^{-1}, & \hbox{$x\in K$;} \\
    0, & \hbox{$x\notin K$.} \\
\end{array}%
\right.
$$
satisfies equation (\ref{l1.1.}), and hence, by Lemma 1 the random
variables $S$ and $D$ are independent.

\textbf{Remark 1.} Let $X$ be a discrete field of characteristic
$2.$ Let $\xi$ and $\eta$ be independent random variables with
values in $X$  and distributions   $\mu$ and $\nu$. If Если the
random variables $S=\xi+\eta$ and $D=(\xi-\eta)^2$ are independent,
then $\mu$ and $\nu$ are degenerate distributions. Indeed, in a
field of characteristic  $2$ the equality
$(\xi-\eta)^2=(\xi+\eta)^2$ holds, i.e. $D=S^2.$ But as it easily
seen, in a filed of characteristic 2 if  $\mu_S$ is a nondegenerate
distribution, then $\mu_{S^2}$ is also a nondegenerate distribution.
Hence, in a filed of characteristic 2 independence of $S$ and $S^2$
implies that $S$ has a degenerate distribution (compare below with
Remark 2). But then $\mu$ and $\nu$ are also degenerate
distributions.

We prove now a theorem on characterization of idempotent
distributions for a discrete field  of characteristic $0$.

\textbf{Theorem 2.} \textit{Let $X$ be a countable discrete field of
characteristic $0$. Let $\xi$ and $\eta$ be independent identically
distributed random variables with values in $X$ and distribution
$\mu$. If the random variables
  $S=\xi+\eta$ and $D=(\xi-\eta)^2$ are independent, then
    $\mu$ is a degenerate distribution.}

\textbf{Proof.} Replacing if it is necessary the random variables
$\xi$ and $\eta$ by new independent random variables $\xi+x$ and
$\eta+x$, we can assume from the beginning that $\mu(0)>0$. Then
Lemma 1 implies that the function  $\mu(x)$ satisfies equation
(\ref{l1.1.}). Since the characteristic of the field  $X$ is 0, $X$
is an Abelian group  with unique division by  2. Then by Lemma 2 the
set $K=\{x\in X: \mu(x)>0\}$ is a subgroup of $X$.

Assume that $\mu(x_0)> 0$ at a point $x_0 \in X$, $x_0\ne 0.$
Consider a subgroup of  $X$ of the form
$$G=\left\{{m x_0\over 2^n}:\ m=0,\pm 1, \pm 2,\dots, \ n=0,1,\dots\right\}.$$
Taking into account that $K$ is a subgroup, equation (\ref{l1.1.})
implies that $\mu(x)>0$ for all $x\in G$. Consider the restriction
of equation (\ref{l1.1.}) to the subgroup $L$.  Put $\varphi(x)=\log
\mu(x)$, $x\in G$. It follows from (\ref{l1.1.}) that
\begin{equation}\label{t2.8}
    2\varphi(0)+\varphi(u+v)+\varphi(u-v)=2 \varphi(u)+\varphi(v)+\varphi(-v), \quad u,v\in G.
\end{equation}
Let $h$ be an arbitrary element of the group $G$. Substituting in
(\ref{t2.8}) $u$ for $u+h$ and $v$ for $v+h$, we obtain
\begin{equation}\label{t2.9}
    2\varphi(0)+\varphi(u+v+2h)+\varphi(u-v)=2 \varphi(u+h)+\varphi(v+h)+\varphi(-v-h),
    \quad u,v,h\in G.
\end{equation}
Subtracting   (\ref{t2.8}) from  (\ref{t2.9})  we find
\begin{equation}\label{t2.10}
    \Delta_{2h}\varphi(u+v)=2 \Delta_h\varphi(u)+\Delta_h\varphi(v)+\Delta_{-h}\varphi(-v),
    \quad u,v,h\in G.
\end{equation}
Putting in (\ref{t2.10}) $u=0$, we get
\begin{equation}\label{t2.11}
    \Delta_{2h}\varphi(v)=2 \Delta_h\varphi(0)+\Delta_h\varphi(v)+\Delta_{-h}\varphi(-v),
    \quad v,h\in G.
\end{equation}
Subtracting  (\ref{t2.11}) from (\ref{t2.10}), we obtain
\begin{equation}\label{t2.12}
    \Delta_{2h}\Delta_{u}\varphi(v)=2 \Delta_h\Delta_{u} \varphi(0),
    \quad u,v,h\in G.
\end{equation}
Let $k$ be an  arbitrary element of the group $G$. Substituting in
(\ref{t2.12})   $v$ for $v+k$ and subtracting from the obtained
equation equation (\ref{t2.12}), we find
\begin{equation}\label{t2.13}
    \Delta_{2h}\Delta_{u}\Delta_k\varphi(v)=0, \quad u,v,h,k\in G.
\end{equation}
Since   $u$, $h$ and $k$ are arbitrary elements of $G$, it follows
from (\ref{t2.13}) that the function $\varphi(x)$ satisfies the
equation
\begin{equation}\label{t2.14}
    \Delta_{u}^3 \varphi(v)=0, \quad u, v\in G.
\end{equation}
Let $x$ be an  arbitrary element of the group $G$.  Then $x=rx_0$,
where $r={m\over 2^n}$. As it easily follows from (\ref{t2.14})
 $\varphi(x)=ar^2+br+c$, where  $a, b, c$
are some constants. Hence, $\mu(x)=e^{ar^2+br+c}$. But this is
impossible because
$$
\sum_{x\in G}\mu(x)\le 1.
$$
Thus, there not exists a point ${x}_0 \in X$, ${x}_0\neq 0$, such
that $\mu({x}_0)> 0$, and hence $\mu=E_0.$

\textbf{Remark 2.} Let $X$ be a countable discrete field of
characteristic $p$, where $p\ne 2$. Let $\xi$ and $\eta$ be
independent   random variables with values in $X$ and distributions
 $\mu$ and $\nu$.  Assume that the random variables
$S=\xi+\eta$ and $D=(\xi-\eta)^2$ are independent. Generally
speaking this not implies that both distributions  $\mu$ and $\nu$
are idempotent. There is a corresponding example.

Since $p\ne 2,$ we have $e\ne -e$. Let $\xi$ and $\eta$ be
independent   random variables with values in $X$ and distributions
 $\mu$ and $\nu$ such that  $\mu={1\over
2}(E_{-e}+E_{e}),$ $\nu=E_0$. Then $S=\xi$ и $D=\xi^2$. Since
$\mu_D$ is a degenerate distribution, the random variables $S$ and
$D$ are independent.

\section{Characterization of idempotent distributions on the field of $p$-adic numbers}
\label{sec:2}

Let $p$ be a fixed prime number. Denote by $\mathbf{Q}_p $ the filed
of $p$-adic numbers. We need some simple properties of the field of
$\mathbf{Q}_p $. As a set $\mathbf{Q}_p $ coincides with a set of
doubly infinite sequences of positive integers
 $$x=(\dots,x_{-n}, x_{-n+1},\dots, x_{-1}, x_0,
x_1,\dots,x_n, x_{n+1},\dots), \ x_n \in\{0, 1,\dots, p-1\}$$, being
$x_n=0$ for $n < n_0$, where the number $n_0$ depends on $x$. We
correspond to each element $x \in \mathbf{Q}_p $ a  series
 $\sum\limits_{k=-\infty}^{\infty} x_k p^k.$ Addition and multiplication
 of   series  are defined in a natural way and define the
 operations
 of addition and multiplication on $\mathbf{Q}_p $.
With respect to these operations  $\mathbf{Q}_p $ is a field. Denote
by $\mathbf{Z}_p $ a subset of  $\mathbf{Q}_p $ consisting of all $x
\in \mathbf{Q}_p $ such that $x_n=0$ for $n < 0$. The set
$\mathbf{Z}_p $ is a ring and is called the ring of $p$-adic
integers.   Denote by $\mathbf{Z}_p^\times $ a subset of
$\mathbf{Q}_p $ of the form $\mathbf{Z}_p^\times =\{x \in
\mathbf{Q}_p : x_n =0$ for $n < 0$, $x_0 \ne 0 \}$. The subset
$\mathbf{Z}_p^\times $ coincides with the group of invertible
elements of the ring $\mathbf{Z}_p $. Moreover each element $x\in
\mathbf{Q}_p $ is uniquely represented in the form
 $x=p^lc$, where $l$ is an integer, and $c\in \mathbf{Z}_p^\times $.
 Denote by $e$ the unity of the field $\mathbf{Q}_p $. One can define a norm
$|x|_p$ on $\mathbf{Q}_p $ by the following way. If $x\in
\mathbf{Q}_p $ и $x=p^lc$, where $l$ is an integer, and  $c\in
\mathbf{Z}_p^\times $, we put $|x|_p=p^{-l}$, $|0|_p=0$. The norm
$|x|_p$ satisfies the conditions: $|x+y|_p\le \max(|x|_p, |y|_p)$,
$|xy|_p=|x|_p|y|_p$ and defines a topology on $\mathbf{Q}_p $. The
field $\mathbf{Q}_p $ with respect to this topology is locally
compact, noncompact and totally disconnected, and the ring
$\mathbf{Z}_p $ is  compact. Choose a Haar measure $m_{\mathbf{Q}_p
}$ on $\mathbf{Q}_p $ in such a way that $m_{\mathbf{Q}_p
}(\mathbf{Z}_p )=1.$ Then $m_{\mathbf{Q}_p }(p^k\mathbf{Z}_p
)=p^{-k}.$ We will also assume that $m_{\mathbf{Q}_p
^2}=m_{\mathbf{Q}_p }\times m_{\mathbf{Q}_p }$.

To prove the main theorem of this section we need some lemmas.

\textbf{Lemma 3.} \textit{Consider the field $\mathbf{Q}_p $. Then
on the set
 $\mathbf{Q}_p ^{[2]}$ there exists a continuous function  $\mathfrak{s}(x)$
satisfying the equation
\begin{equation}\label{eq1}
 \mathfrak{s}^2(x)=x, \quad x\in\mathbf{Q}_p ^{[2]}.
\end{equation}
The set  $\mathbf{Q}_p^{[2]} \backslash\{0\}$  is a disjoin union of
balls in each of which the function $\mathfrak{s}(x)$ is expressed
by a convergent power series.}

 \textbf{Proof.}  The representation
$$\mathbf{Q}_p =\{0\}\cup\bigcup_{l=-\infty}^{\infty}p^{l}\mathbf{Z}_p^\times
$$
implies that
\begin{equation}\label{eq2}
\mathbf{Q}_p
^{[2]}=\{0\}\cup\bigcup_{l=-\infty}^{\infty}p^{2l}(\mathbf{Z}_p^\times
)^{[2]}.
\end{equation}

Assume first that  $p>2$.  It is well known that in the
multiplicative group  $\mathbf{Z}_p^\times $ there exists an element
$\varepsilon$ of order $p-1$, and elements  $0, \varepsilon,
\varepsilon^2, \dots, \varepsilon^{p-1}=e$ form a complete set of
coset representatives of the subgroup $p\mathbf{Z}_p $ in  the group
$\mathbf{Z}_p $. Put
$A_k=\varepsilon^{k}+p\mathbf{Z}_p=\left\{x:|x-\varepsilon^{k}|_p\le
{1\over p}\right\}$, $k=1, 2,\dots, p-1$. Then $\mathbf{Z}_p^\times
=\bigcup\limits_{k=1}^{p-1}A_{k}.$ First define the function
$\mathfrak{s}(x)$ on the set $(\mathbf{Z}_p^\times )^{[2]}$. Since
$(e+p\mathbf{Z}_p )^{[2]}=e+p\mathbf{Z}_p $, we have
$A_k^{[2]}=A_{2k}$, if $k=1,2,\dots, {p-1\over 2},$ and
$A_k^{[2]}=A_{2k - p+1}$, if $k={p+1\over 2}, {p+3\over 2}, \dots,
{p-1}.$ Note that $(\mathbf{Z}_p^\times
)^{[2]}=\bigcup\limits_{k=1}^{p-1\over2}A_{2k},$ and  define the
function $\mathfrak{s}(x)$ on each coset $A_{2k}$. Take $x\in
A_{2k}$. Then the equation $x=t^2$ has two solutions $t_1\in A_k$
and $-t_1\in A_{k+{p-1\over 2}}$. These solutions belong to
different cosets. The coset
 $A_k$
is a compact set, and the function $g(x)=x^2$ is continuous on $A_k$
and it is one-to-one mapping of the set  $A_k$ on $A_{2k}$. This
implies that the  inverse to  $g(x)$ mapping
$\mathfrak{s}_k:A_{2k}\mapsto A_k$ is also continuous,  and hence is
a homeomorphism between $A_{2k}$ and $A_k$.  Put
$\mathfrak{s}(x)=\mathfrak{s}_k(x)$, if $x\in A_{2k}$, $k=1,
2,\dots, {p-1\over 2}$. Since $A_{2k}$ is an open set in
$\mathbf{Q}_p $, the function $\mathfrak{s}(x)$ is continuous and
satisfies equation (\ref{eq1}) on $(\mathbf{Z}_p^\times )^{[2]}$.
Taking into account (\ref{eq2}), put
$$ \mathfrak{s}(x)=
\left\{%
\begin{array}{ll}
    p^l\mathfrak{s}(c), & \hbox{$x=p^{2l}c, c\in (\mathbf{Z}_p^\times )^{[2]}$;} \\
    0, & \hbox{x=0.} \\
\end{array}%
\right.
$$
It is not difficult to verify that the constructed function
$\mathfrak{s}(x)$ in each ball $A_{2k}$ при $k=1, 2,\dots, {p-1\over
2}$ is expressed by a convergent power series
$$
\mathfrak{s}(x)=\varepsilon^k\left(e+{\varepsilon^{-2k}\over 2}
(x-\varepsilon^{2k})+\sum_{n=2}^\infty(-1)^{n-1}{(2n-3)!!\over
2n!!}\varepsilon^{-2kn}(x-\varepsilon^{2k})^n\right).
$$
This implies that in each ball  $p^{2l}A_{2k}$, $l= \pm 1, \pm 2,
\dots$ the function $\mathfrak{s}(x)$ is also expressed by a
convergent power series. It is obvious that the set
$\mathbf{Q}_p^{[2]} \backslash\{0\}$ is   a disjoin union of balls
$p^{2l}A_{2k}$, $k=1, 2,\dots, {p-1\over 2}$, $l=0, \pm 1, \pm 2,
\dots$, and  $\mathfrak{s}(x)$ is the required function.

If $p=2$, the reasoning is changed slightly. Put
$B_k=ke+4\mathbf{Z}_2$, $k=0, 1, 2, 3$. Then the elements  $0, e,
2e, 3e$ form a complete set of coset representatives of the subgroup
$4\mathbf{Z}_2$ in the group $\mathbf{Z}_2$. We have
$(\mathbf{Z}_2^\times)^{[2]}=B_1^{[2]}=B_3^{[2]}=e+8\mathbf{Z}_2$.
Take $x\in e+8\mathbf{Z}_2$. Then the equation $x=t^2$ has two
solutions $t_1\in B_1$ and $-t_1\in B_3$. These solutions belong to
different cosets. The rest part of the proof is similar to the case
when $p>2$.

\textbf{Lemma 4.} \textit{Consider the field $\mathbf{Q}_p $, where
$p>2$. Take $(x_0, y_0)\in \mathbf{Q}_p ^2$ such that
$|x_0-y_0|_p=p^{-l}.$ Then for  $k\ge l+1$  the following equality
\begin{equation}\label{eq16}
T\{(x_0, y_0)+(p^{k}\mathbf{Z}_p )^2\}=(x_0+y_0,
(x_0-y_0)^2)+(p^{k}\mathbf{Z}_p )\times(p^{k+l}\mathbf{Z}_p )
\end{equation}
holds}.

\textbf{Proof.} Since $|x_0-y_0|_p=p^{-l}$, we have $x_0-y_0=p^lc$,
where $c\in \mathbf{Z}_p^\times $. Note that on the one hand,
\begin{equation}\label{eq15}
T\{(x_0, y_0)+(p^{k}\mathbf{Z}_p )^2\}=T\{(x_0+p^{k}x, y_0+p^{k}y):
x, y\in \mathbf{Z}_p \}=$$$$=\{(x_0+y_0+p^{k}(x+y), (x_0-y_0)^2+
2p^{k}(x_0-y_0)(x-y)+p^{2k}(x-y)^2): x, y\in \mathbf{Z}_p \}=$$$$=
\{(x_0+y_0+p^{k}s, (x_0-y_0)^2+ 2p^{k+l}ct+p^{2k}t^2): s, t\in
\mathbf{Z}_p \}
\end{equation}
holds true for any $k$. On the other hand,
\begin{equation}\label{eq12}
\{2ct+p^{k-l}t^2: t\in \mathbf{Z}_p \}=\mathbf{Z}_p
\end{equation}
holds true for $k\ge l+1$.

Indeed, note that $(e+p^m\mathbf{Z}_p )^{[2]}=e+p^m\mathbf{Z}_p $ is
fulfilled for any
 $m\ge 1$. This implies that
$(c+p^m\mathbf{Z}_p )^{[2]}=c^2+p^m\mathbf{Z}_p $ for all $c\in
\mathbf{Z}_p^\times $, i.e. $\{c^2+2cp^mt+p^{2m}t^2:t\in
\mathbf{Z}_p \}=c^2+p^m\mathbf{Z}_p ,$ and hence,
$\{2ct+p^{m}t^2:t\in \mathbf{Z}_p \}=\mathbf{Z}_p .$ For $k\ge l+1$
this equality implies (\ref{eq12}). Taking into account
(\ref{eq12}), (\ref{eq16}) follows from (\ref{eq15}).

\textbf{Lemma 5.} \textit{Consider the field $\mathbf{Q}_p $, where
$p>2$. Let a function $\mathfrak{s}$ be as constructed in the proof
of Lemma $3$. Consider the mappings
 $S_j$ from $\mathbf{Q}_p \times \mathbf{Q}_p ^{[2]}$ to $\mathbf{Q}_p ^2$ of the form
$$
S_1(u, v)=\left({u+\mathfrak{s}(v)\over 2}, {u-\mathfrak{s}(v)\over
2}\right), \quad S_2(u, v)=\left({u-\mathfrak{s}(v)\over 2},
{u+\mathfrak{s}(v)\over 2}\right).
$$
Let $(u_0, v_0)\in\mathbf{Q}_p \times \mathbf{Q}_p ^{[2]}$ and
$|\mathfrak{s}(v_0)|_p=p^{-l}.$ Put $E_k= \{(u_0,
v_0)+(p^{k}\mathbf{Z}_p )\times(p^{k+l}\mathbf{Z}_p )\}$. Then for
$k\ge l+1$ the following statements are valid:
$$E_k\subset \mathbf{Q}_p \times \mathbf{Q}_p ^{[2]}, \leqno (i)$$
$$S_1(E_k)\cap S_2(E_k)=\emptyset, \leqno (ii)$$
$$
\int\limits_{S_j(E_k)}\Phi_j(x, y)dm_{\mathbf{Q}_p
^2}(x,y)=\int\limits_{E_k}\Phi_j(S_j(u,
v))|\mathfrak{s}(v)|_p^{-1}dm_{\mathbf{Q}_p ^2}(u, v), \ j=1, 2,
\leqno (iii)
$$
for any continuous function $\Phi_j(x, y)$ on $S_j(E_k)$.}

\textbf{Proof.} Note that $|v_0|_p=p^{-2l}$. It follows from the
proof of Lemma 3 that if
 $w_0\in \mathbf{Q}_p ^{[2]}$ and
$|w_0|_p=p^{-2l}$, then $w_0+w\in \mathbf{Q}_p ^{[2]}$ for $w\in
p^{2l+1}\mathbf{Z}_p $.  Since    $k\ge l+1$, from what has been
said it follows that $(i)$ is fulfilled.

To prove  $(ii)$ assume that $S_1(E_k)\cap S_2(E_k)\not=\emptyset$.
Then as easily seen, there exist elements
 $v_1, v_2\in \mathbf{Z}_p $ such that
\begin{equation}\label{e201}
\mathfrak{s}(v_0+p^{k+l}v_1)+\mathfrak{s}(v_0+p^{k+l}v_2)\in
p^k\mathbf{Z}_p .
\end{equation}
Since $v_0\in \mathbf{Q}_p ^{[2]}$ and $|v_0|_p=p^{-2l}$, we have
$v_0=p^{2l}c,$ where $c\in (\mathbf{Z}_p^\times )^{[2]}$ and
$v_0+p^{k+l}v_i=p^{2l}(c+p^{k-l}v_i)$, $i=1,2$. This implies that
$\mathfrak{s}(v_0+p^{k+l}v_1)+\mathfrak{s}(v_0+p^{k+l}v_2)=
p^l(\mathfrak{s}(c+p^{k-l}v_1)+\mathfrak{s}(c+p^{k-l}v_2)).$ It
follows from the definition of the function $\mathfrak{s}$ that the
elements $\mathfrak{s}(c+p^{k-l}v_i)$, $i=1, 2,$ are at the same
coset of the subgroup $p\mathbf{Z}_p$ in the group $\mathbf{Z}_p $,
and hence,
$\mathfrak{s}(c+p^{k-l}v_1)+\mathfrak{s}(c+p^{k-l}v_2)\not\in
p\mathbf{Z}_p.$ Therefore
$\mathfrak{s}(v_0+p^{k+l}v_1)+\mathfrak{s}(v_0+p^{k+l}v_2)\not\in
p^{l+1}\mathbf{Z}_p ,$ But this contradicts (\ref{e201}) for $k\ge
l+1$. So, we proved $(ii)$.

Let us prove  $(iii)$. We will prove that equality $(iii)$ holds
true for $S_1$. For $S_2$ the reasoning is similar. Put $(x_0,
y_0)=S_1(u_0, v_0)=\left({u_0+\mathfrak{s}(v_0)\over 2},
{u_0-\mathfrak{s}(v_0)\over 2}\right)$ and check that
\begin{equation}\label{e266}
S_1(E_k)=(x_0, y_0)+(p^{k}\mathbf{Z}_p )^2.
\end{equation}
We have $S_1(u_0+u, v_0+v)=\left({u_0+u+\mathfrak{s}(v_0+v)\over 2},
{u_0+u-\mathfrak{s}(v_0+v)\over 2}\right)$. It is easily seen that
\begin{equation}\label{e263}
\left|{u_0+u+\mathfrak{s}(v_0+v)\over 2}-
{u_0+\mathfrak{s}(v_0)\over 2}\right|_p\le\max{\{|u|_p,
|\mathfrak{s}(v_0+v)-\mathfrak{s}(v_0)|_p\}}.
\end{equation}
Under the condition of the theorem
\begin{equation}\label{e264}
|u|_p\le p^{-k}.
\end{equation}
Since $\mathfrak{s}^2(x)=x$, we have
\begin{equation}\label{e261}
|\mathfrak{s}(v_0+v)-\mathfrak{s}(v_0)|_p={|v|_p\over|\mathfrak{s}(v_0+v)+\mathfrak{s}(v_0)|_p}.
\end{equation}
Note that $v_0+v=p^{2l}c+p^{k+l}t$ for some $t\in \mathbf{Z}_p $.
This implies that $\mathfrak{s}(v_0+v)=p^l\mathfrak{s}(c+p^{k-l}t)$
and moreover, $\mathfrak{s}(v_0)=p^l\mathfrak{s}(c)$. Inasmuch as
the points $\mathfrak{s}(c+p^{k-l}t)$ and  $\mathfrak{s}(c)$ are at
the same coset of the subgroup    $p\mathbf{Z}_p $ in the group
$\mathbf{Z}_p $, so is
$|\mathfrak{s}(c+p^{k-l}t)+\mathfrak{s}(c)|_p=1$, and hence,
 $|\mathfrak{s}(v_0+v)+\mathfrak{s}(v_0)|_p=p^{-l}$.
If it is remembered that $|v|_p\le p^{-k-l}$, we get from
(\ref{e261}) that
\begin{equation}\label{e262}
|\mathfrak{s}(v_0+v)-\mathfrak{s}(v_0)|_p\le p^{-k}.
\end{equation}
Taking into account (\ref{e264}) and (\ref{e262}), it follows from
(\ref{e263}) that inequality
$$
\left|{u_0+u+\mathfrak{s}(v_0+v)\over 2}-
{u_0+\mathfrak{s}(v_0)\over 2}\right|_p\le p^{-k}
$$
holds, and hence
\begin{equation}\label{e265}
S_1(E_k)\subset (x_0, y_0)+(p^{k}\mathbf{Z}_p)^2.
\end{equation}

We note that if $T(a, b)=T(a', b')$, then either $(a, b)=(a', b')$,
or $(a, b)=(b', a')$. Since
$|x_0-y_0|_p=|\mathfrak{s}(v_0)|_p=p^{-l},$ and $k\ge l+1$, the
restriction of the mapping $T$ to the set $(x_0,
y_0)+(p^{k}\mathbf{Z}_p)^2$ is injective.  Taking this into account,
(\ref{e266}) follows from Lemma 4 and (\ref{e265}) . Moreover it
follows from has been said the mappings  $T$ and $S_1$ are inverse
homeomorphisms of the sets $S_1(E_k)$ and $E_k$.

Observe that the set $E_k$ is a product of the balls
$E_k=\left\{u:|u-u_0|_p\le {1\over p^k}\right\}\times
\left\{v:|v-v_0|_p\le {1\over p^{k+l}}\right\}$. Put $x(u,
v)={u+\mathfrak{s}(v)\over 2},$ $y(u, v)={u-\mathfrak{s}(v)\over
2}$, $(u, v)\in E_k$. The mapping $S_1$ is a homeomorphism of the
open compacts $E_k$  and $S_1(E_k)$, and the functions $x(u, v)$ and
$y(u, v)$ by Lemma 3 on $E_k$ are expressed by convergent power
series in $u$ and $v$.  We have ${\mathbf{d}x\over \mathbf{d}
u}={\mathbf{d}y\over \mathbf{d} u}={1\over 2},$ ${\mathbf{d}x\over
\mathbf{d} v}={1\over 4\mathfrak{s}(v)},$ ${\mathbf{d}y\over
\mathbf{d} v}=-{1\over 4\mathfrak{s}(v)}.$ It follows from this that
\begin{equation}\label{eq28.02}
\left|\rm det\left(%
\begin{array}{cc}
  {\mathbf{d}x\over \mathbf{d}
u} & {\mathbf{d}x\over \mathbf{d} v} \\
  {\mathbf{d}y\over \mathbf{d} u} & {\mathbf{d}y\over \mathbf{d} v} \\
\end{array}%
\right)\right|_p=\left|\rm det\left(%
\begin{array}{cc}
  {1\over 2} & {1\over 4\mathfrak{s}(v)} \\
 {1\over 2}  & -{1\over 4\mathfrak{s}(v)} \\
\end{array}%
\right)\right|_p=|\mathfrak{s}(v)|_p^{-1}\ne 0, \quad (u, v)\in E_k.
\end{equation}
It is obvious that $(iii)$ follows from the change of variables
formula in integrals (\cite[\S 4]{VVZ}) and  (\ref{eq28.02}).

\textbf{Lemma 6.}  \textit{Consider the field $\mathbf{Q}_p $, where
$p>2$. Let $\xi$ and $\eta$ be independent identically distributed
random variables with values in $\mathbf{Q}_p $   and distribution
 $\mu$. Assume that $\mu$ has the density  $\rho$
with respect to  $m_{\mathbf{Q}_p }$ such that $\rho$ is continuous
and $\rho(0)>0$. The random variables
  $S=\xi+\eta$ and $D=(\xi-\eta)^2$ are independent if and only if the density
 $\rho$ satisfies the equation}
\begin{equation}\label{eq22}
\rho^2(u)\rho(v)\rho(-v)=\rho^2(0)\rho(u+v)\rho(u-v), \quad u,v\in
    \mathbf{Q}_p .
\end{equation}

\textbf{Proof.}  Inasmuch as $\mu_{(S,D)}=T(\mu_{(\xi,\eta)})$  and
the distribution $\mu_{(\xi,\eta)}$ is absolutely continuous with
respect to $m_{\mathbf{Q}_p ^2}$, so is $\mu_{(\xi,\eta)}\{(t, t):
t\in \mathbf{Q}_p \}=0$. Therefore the distribution  $\mu_{(S,D)}$
concentrated at the set
 $\mathbf{Q}_p \times (\mathbf{Q}_p ^{[2]}\backslash \{0\})$.
Let   the mappings $S_j$ and the sets $E_k$ be the same as in Lemma
5. Take $(u_0, v_0)\in \mathbf{Q}_p \times (\mathbf{Q}_p
^{[2]}\backslash \{0\})$ and represent the element
$\mathfrak{s}(v_0)$ in the form $\mathfrak{s}(v_0)=p^lc$, where
$c\in \mathbf{Z}_p^\times $. By Lemma $(i)$ and $(ii)$ are fulfilled
for $k\ge l+1$. We have
\begin{equation}\label{eq4}
\mu_{(S,D)}\{E_k\}=T(\mu_{(\xi,\eta)})\{E_k\}=\mu_{(\xi,\eta)}
\{T^{-1}(E_k)\}=\int\limits_{T^{-1}(E_k)}\rho(x)\rho(y)dm_{\mathbf{Q}_p
^2}(x,y)=$$$$ \int\limits_{S_1(E_k)}\rho(x)\rho(y)dm_{\mathbf{Q}_p
^2}(x,y)+\int\limits_{S_2(E_k)}\rho(x)\rho(y)dm_{\mathbf{Q}_p
^2}(x,y).
\end{equation}
Taking into account equalities  $(iii)$ of Lemma 5, we transform
integrals in the right-hand side of equality (\ref{eq4}) and obtain
$$
\int\limits_{S_1(E_k)}\rho(x)\rho(y)dm_{\mathbf{Q}_p
^2}(x,y)=\int\limits_{E_k}\rho\left({u+\mathfrak{s}(v)\over
2}\right) \rho\left({u-\mathfrak{s}(v)\over
2}\right)|\mathfrak{s}(v)|_p^{-1}dm_{\mathbf{Q}_p ^2}(u, v),
$$
$$
\int\limits_{S_2(E_k)}\rho(x)\rho(y)dm_{\mathbf{Q}_p
^2}(x,y)=\int\limits_{E_k}\rho\left({u-\mathfrak{s}(v)\over
2}\right) \rho\left({u+\mathfrak{s}(v)\over
2}\right)|\mathfrak{s}(v)|_p^{-1}dm_{\mathbf{Q}_p ^2}(u, v).
$$
Then (\ref{eq4}) implies that
$$
\mu_{(S,D)}\{E_k\}=2\int\limits_{E_k}\rho\left({u+\mathfrak{s}(v)\over
2}\right) \rho\left({u-\mathfrak{s}(v)\over
2}\right)|\mathfrak{s}(v)|_p^{-1}dm_{\mathbf{Q}_p ^2}(u, v).
$$
This equality means that the distribution  $\mu_{(S,D)}$ has a
density with respect to $m_{\mathbf{Q}_p ^2}$, and this density is
equal to
\begin{equation}\label{eq20}
2\rho\left({u+\mathfrak{s}(v)\over
2}\right)\rho\left({u-\mathfrak{s}(v)\over 2}\right)
|\mathfrak{s}(v)|_p^{-1}, \quad u\in\mathbf{Q}_p , \
v\in\mathbf{Q}_p ^{[2]}\backslash \{0\}.
\end{equation}

Note that when we got representation (\ref{eq20}) for the density of
distribution $\mu_{(S,D)}$, we did not use the independence of the
random variables  $S$ and $D$.

By the condition of the lemma the random variables $S$ and $D$ are
independent. Therefore there exist integrable with respect to
$m_{\mathbf{Q}_p }$ functions $r_j$ on $\mathbf{Q}_p $ such that the
equality
\begin{equation}\label{eq5}
r_1(u)r_2(v)=2\rho\left({u+\mathfrak{s}(v)\over
2}\right)\rho\left({u-\mathfrak{s}(v)\over
2}\right)|\mathfrak{s}(v)|_p^{-1}
\end{equation}
holds true almost everywhere on $\mathbf{Q}_p \times (\mathbf{Q}_p
^{[2]}\backslash \{0\})$ with respect to $m_{\mathbf{Q}_p ^2}$.
Since the function in the right-hand side of equality (\ref{eq5}) i
continuous, we can assume without loss of generality that the
functions  $r_j$ are also continuous, and equality  (\ref{eq5})
holds true everywhere on $\mathbf{Q}_p \times (\mathbf{Q}_p
^{[2]}\backslash \{0\})$. Since $\rho(0)>0$, it is easily seen that
$r_1(0)>0$. Put $v=t^2$, $t\neq 0$. It follows from  (\ref{eq5})
that
\begin{equation}\label{eq6}
r_2(t^2)=2r_1^{-1}(0)\rho\left({\mathfrak{s}(t^2)\over
2}\right)\rho\left(-{\mathfrak{s}(t^2)\over 2}\right)
|\mathfrak{s}(t^2)|_p^{-1}, \  t \in\mathbf{Q}_p , \quad t\neq 0.
\end{equation}
Note that  (\ref{eq5}) and (\ref{eq6}) imply the equality
\begin{equation}\label{eq7}
r_1(u) \rho\left({\mathfrak{s}(t^2)\over
2}\right)\rho\left(-{\mathfrak{s}(t^2)\over 2}\right)=
r_1(0)\rho\left({u+\mathfrak{s}(t^2)\over
2}\right)\rho\left({u-\mathfrak{s}(t^2)\over 2}\right), $$ $$\quad
(u, t)\in\mathbf{Q}_p ^2, \ t\neq 0.
\end{equation}
It follows from the continuity of  $\rho$ and $r_1$ that equality
(\ref{eq7}) holds true for all $u, t\in\mathbf{Q}_p $. Put in
(\ref{eq7}) $t=0$. We deduce from the last equality that
\begin{equation}\label{eq8}
r_1(u)={r_1(0)\over \rho^2(0)}\rho^2\left({u\over 2}\right), \quad
u\in\mathbf{Q}_p .
\end{equation}
Substituting  (\ref{eq8}) into (\ref{eq7}), we find that
\begin{equation}\label{eq9}
\rho^2\left({u\over 2}\right)\rho\left({\mathfrak{s}(t^2)\over
2}\right)\rho\left(-{\mathfrak{s}(t^2)\over 2}\right)=
\rho^2(0)\rho\left({u+\mathfrak{s}(t^2)\over
2}\right)\rho\left({u-\mathfrak{s}(t^2)\over 2}\right), $$ $$\quad
u, t\in\mathbf{Q}_p .
\end{equation}
Since either     $\mathfrak{s}(t^2)=t$, or  $\mathfrak{s}(t^2)=-t$,
the equalities
\begin{equation}\label{eq10}
\rho\left({\mathfrak{s}(t^2)\over
2}\right)\rho\left(-{\mathfrak{s}(t^2)\over 2}\right)=
\rho\left({t\over 2}\right)\rho\left(-{t\over 2}\right), \quad   t
\in\mathbf{Q}_p
\end{equation}
and
\begin{equation}\label{eq11}
\rho\left({u+\mathfrak{s}(t^2)\over
2}\right)\rho\left({u-\mathfrak{s}(t^2)\over 2}\right)=
\rho\left({u+t\over 2}\right)\rho\left({u-t\over 2}\right), \quad u,
t\in\mathbf{Q}_p
\end{equation}
are fulfilled. Substituting  (\ref{eq10}) and (\ref{eq11}) into
(\ref{eq9}), we get that the density  $\rho$ satisfies equation
(\ref{eq22}). The necessity is proved.

Let us prove the sufficiency. Taking into account  (\ref{eq20}) and
(\ref{eq11})  we have the following representation for the density
$\varrho$ of the distribution $\mu_{(S,D)}$:
\begin{equation} \label{eq21} \varrho(u, v)=
\left\{%
\begin{array}{ll}
    2\rho\left({u+t\over 2}\right)\rho\left({u-t\over
2}\right)|\mathfrak{s}(t^2)|_p^{-1}, & \hbox{$u\in
\mathbf{Q}_p , \ v=t^2, \ t\ne 0$;} \\
    0, & \hbox{$u\in \mathbf{Q}_p , \ v\notin(\mathbf{Q}_p ^{[2]}\backslash
\{0\})$.} \\
\end{array}%
\right.
\end{equation}
If a density  $\rho$  satisfies equation  (\ref{eq22}), it is easily
seen that the density
 $\varrho(u, v)$ is represented aa a product of a function of $u$ and a function of  $v$.
 This implies the independence of  $S$ and $D$.

Now we can prove the main theorem of this section.

\textbf{Theorem 3.} \textit{Consider the field $\mathbf{Q}_p $,
where $p>2$. Let $\xi$ and $\eta$ be independent identically
distributed random variables with values in $\mathbf{Q}_p $   and
distribution
 $\mu$. Assume that $\mu$ has the density  $\rho$
with respect to  $m_{\mathbf{Q}_p }$ such that $\rho$ is continuous.
In order that the random variables
  $S=\xi+\eta$ and $D=(\xi-\eta)^2$ be independent it is necessary
  and sufficient that  $\mu\in I(\mathbf{Q}_p )$.}

\textbf{Proof.} Necessity. We reason as in the proof of Theorem 1.
It is obvious that replacing if it is necessary the random variables
$\xi$ and $\eta$ by new independent random variables $\xi+x$ and
$\eta+x$, we can assume from the beginning that $\rho(0)>0$.
Applying  Lemma 6 we find  that the density $\rho$ satisfies
equation (\ref{eq22}). Since the field   $\mathbf{Q}_p $ is an
Abelian group  with unique division by  2, by Lemma 2 the set
$K=\{x\in \mathbf{Q}_p : \rho(x)>0\}$ is a subgroup of $\mathbf{Q}_p
$. Obviously, $K$ is an open subgroup, and hence a closed one.

Assume first that $K\ne \mathbf{Q}_p $. Since $K\ne\{0\}$, the
subgroup
 $K$ is of the form $K=p^m\mathbf{Z}_p $ for some integer
  $m$, and hence, $K$
is a compact group. Consider the restriction of equation
(\ref{eq22}) to $K$. Put $\varphi(x)=\log \rho(x)$, $x\in K$. It
follows from (\ref{eq22}) that
$$
2 \varphi(u)+\varphi(v)+\varphi(-v)=
2\varphi(0)+\varphi(u+v)+\varphi(u-v), \quad u ,v \in K.
$$
Integrate both sides of this equality by the measure $dm_K(v)$. We
get that $\varphi(u)=\varphi(0)$ for all $u\in K$. Hence,
$\rho(x)=\rho(0)$ for all $x\in K$, but this means that $\mu=m_K$.

If $K=\mathbf{Q}_p $, then we consider the restriction of equation
(\ref{eq22}) to an arbitrary subgroup $G=p^m\mathbf{Z}_p $ of
$\mathbf{Q}_p $. The reasoning  above shows that $\rho(x)=\rho(0)$
for all $x\in G,$ and hence, $\rho(x)=\rho(0)$ for all $x\in
\mathbf{Q}_p ,$ contrary to the integrability of the density $\rho.$
So, the case when $K=\mathbf{Q}_p $, is impossible. The necessity is
proved.

Sufficiency.  Let $K$ be a nonzero compact subgroup of $\mathbf{Q}_p
$. As has been stated above, $K$ is of the form
 $K=p^m\mathbf{Z}_p $ for some integer $m$. Let $\xi$ and $\eta$ be
independent identically distributed random variables with values in
$\mathbf{Q}_p $   and distribution
 $\mu=m_K$.  Since $p>2$, we have ${x\over
2}\in K$ for any element $x$ of the subgroup $K$. It is well known
that this implies that the random variables $\xi+\eta$ and
$\xi-\eta$ are independent (\cite[\S 7]{Fe1})). Hence, the random
variables $S$ and  $D$ are also independent.

\textbf{Remark 3.} Let us discuss the case of the field
$\mathbf{Q}_2 $. A lemma similar to Lemma 4 holds also true for the
field $\mathbf{Q}_2 $. Unlike equality  (\ref{eq16}) for $k\ge l+2$
the equality
\begin{equation}\label{eq30}
T\{(x_0, y_0)+(2^{k}\mathbf{Z}_2)^2\}=(x_0+y_0,
(x_0-y_0)^2)+(2^{k}\mathbf{Z}_2)\times(2^{k+l+1}\mathbf{Z}_2)
\end{equation}
holds true. Respectively, taking (\ref{eq30}) into account, one can
also reformulate Lemma 5. It allows to prove Lemma 6 for the field
 $\mathbf{Q}_2 $. Next reasoning as in the proof of the necessity in Theorem 3, we get
 that if $\xi$ and $\eta$ are
independent identically distributed random variables with values in
$\mathbf{Q}_2 $   and distribution
 $\mu$ such that $\mu$ has the density  $\rho$
with respect to  $m_{\mathbf{Q}_p }$, where the density $\rho$ is
continuous, $\rho(0)>0$ and
 the random variables
  $S=\xi+\eta$ and $D=(\xi-\eta)^2$ are independent, then
$\mu=m_K$, where $K=2^m\mathbf{Z}_2$ for some integer $m$. But in
this case the corresponding density
$$ \rho(x)=
\left\{%
\begin{array}{ll}
    2^m, & \hbox{$x\in 2^m\mathbf{Z}_2$;} \\
    0, & \hbox{$x\notin 2^m\mathbf{Z}_2$.} \\
\end{array}%
\right.
$$
does not satisfies equation (\ref{eq22}). To see this, put in
(\ref{eq22}) $u=v\in 2^{m-1}\mathbf{Z}_2\backslash 2^m\mathbf{Z}_2$.
Then the left-hand side of (\ref{eq22}) is equal to zero, but the
right-hand side does not. Thus, do not exist independent identically
distributed random variables
 $\xi$ and $\eta$ with values
in the field $\mathbf{Q}_2 $ and distribution $\mu$, such that $\mu$
has the density  $\rho$ with respect to a Haar measure
 $m_{\mathbf{Q}_2 }$, where the density $\rho$ is continuous and
 the random variables
  $S=\xi+\eta$ and $D=(\xi-\eta)^2$ are independent.

\section{Comments and unsolved problems}
\label{sec:3}

Let $\xi$ and  $\eta$ be independent  random variables with values
in the field of real numbers $\mathbf{R}$. According to the
Kac-Bernstein theorem the independence of  $\xi+\eta$  and
$\xi-\eta$ implies that the random variables  $\xi$ and $\eta$ are
Gaussian. The similar result holds also  true for the considering in
the article fields. Indeed, let $\xi$ and  $\eta$ be independent
random variables with values either in a countable discrete field or
in the field $\mathbf{Q}_p $ and with distributions  $\mu$ and
$\nu$. Since in the both cases the connected component of zero of
the field does not contain elements of order 2, the independence of
$\xi+\eta$  and $\xi-\eta$ implies that  $\mu$ and $\nu$ are
idempotent distributions (\cite[Theorem 7.10]{Fe1}). In so doing we
do net assume that $\xi$ and  $\eta$ are identically distributed,
and in the case when $\xi$ and $\eta$ take values
 in   the field $\mathbf{Q}_p $ we do not assume that $\mu$ and $\nu$ have continuous densities.

The situation changes if instead of  $\xi+\eta$  and $\xi-\eta$ we
will consider $S=\xi+\eta$ and $D=(\xi-\eta)^2$. As has been proved
in \cite{KLR}, if one consider independent  random variables with
values in the field of real numbers $\mathbf{R}$, then without the
condition of identically distributivity of  $\xi$ and $\eta$ the
independence of  $S$ and $D$, generally speaking, does not imply
that both distributions of random variables
 $\xi$ and $\eta$ are Gaussian. In  \cite{KLR} was constructed an example, when
 one of the random variables has a Gaussian distribution, but the distribution of
 the second random variable is a mixture of two Gaussian distributions.
It follows from Remark 2 that in Theorems 1 and 2  one can not omit
the condition of identically distributivity of independent  random
variables $\xi$ and $\eta$. The problem, if we can not omit the
condition of identically distributivity of independent  random
variables $\xi$ and $\eta$ in Theorem 3 is unsolved.

The problem, if we can prove Theorem 3 without assumption that the
distribution
 $\mu$ has a density with respect to  $m_{\mathbf{Q}_p }$  and this density is continuous,
is also unsolved.

Finally, the theorems proved in the article are analogues for the
considering fields of a particular case of the theorem Geary. It is
interesting to find out if the theorem Geary holds true for the
considering fields in the general case, i.e. when  $n\ge 2$.



\begin{thebibliography}{}
%
%
\bibitem{Fe10} G.M. Feldman, { On the Skitovich-Darmois
theorem on Abelian groups}. Theory Probab. Appl. \textbf{
    37} (1992),  621--631.


\bibitem{Fe3} G.M. Feldman. { A characterization of the Gaussian
distribution on Abelian groups}. Probab. Theory Relat. Fields, {\bf
126},  (2003), 91--102.

\bibitem{Fe18} G.M. Feldman. { On a characterization theorem for locally
compact abelian groups}.   Probab. Theory Relat. Fields,  {\bf 133},
(2005), 345--357.

\bibitem{Fe19}  G.M. Feldman. { On the Heyde theorem for discrete Abelian
groups}.   Studia Math., {\bf 177},   (2006), 67--79.

\bibitem{Fe1} G.M. Feldman.  Functional equations and characterization problems on locally
compact Abelian groups. EMS Tracts in Mathematics, {\bf 5}, Zurich:
European Mathematical Society (EMS), 2008, 268 p.

\bibitem{Fe30}  G.M. Feldman. { The Heyde theorem for locally compact Abelian
groups}, J. of Funct. Analysis  {\bf 258},  (2010), 3977--3987.

\bibitem{JOTP}  G.M. Feldman.  On the Skitovich–Darmois theorem for the group of $p$-adic numbers. J. of Theor. Probability. DOI  10.1007/ s10959-013-0525-9


\bibitem{Ge}  R. C. Geary. The distribution of ``Student's" ratio for non-normal
samples. Suppl. J. R. statist. Soc., London, {\bf 3}, (1936),
178-184.

\bibitem{GraLo}  P.  Graczyk, J.-J. Loeb.  A Bernstein property of measures
on groups and symmetric spaces. Probab. Math. Statist,  \textbf{20},
No 1, (2000), 141--149.

\bibitem{KLR} A. Kagan, R.C. Laha, V. Rohatgi.
Independence of the sum and absolute difference of independent
random variables does not imply their normality. Math. Methods
Statist., {\bf 6}, No. 2,  (1997), 263–265.

\bibitem{KaSa}  T. Kawata, and H. Sakamoto. On the characterisation of
the normal population by the independence of the sample mean and the
sample variance. J. Math. Soc. Japan, {\bf 1}, (1949), 111--115.


\bibitem{Lu}  E. Lukacs. A characterization of the normal
distribution. Ann. Math. Statistics, {\bf 13}, (1942), 91--93.

\bibitem{M} I. P. Mazur.   Skitovich–Darmois theorem for discrete and compact totally disconnected Abelian groups. Ukrainian Mathematical Journal,
{\bf 65}, no 7,  (2013), 1054-1070.

\bibitem{My} M. Myronyuk. The Heyde theorem on $a$-adic solenoids, Colloquium Mathematicum,  {\bf 132}, (2013), 195--210.

\bibitem{NeueScho}   D. Neuenschwander and
R. Schott. {  The Bernstein and Skitovich-Darmois characterization
theorems for Gaussian distributions on groups, symmetric spaces, and
quantum groups}. Exposition. Math., {\bf 15},  (1997), 289--314.

\bibitem{VVZ}  V.S. Vladimirov, I.V. Volovich,  E.I. Zelenov. p-adic analysis and mathematical physics.
Series on Soviet and East European Mathematics. 10. Singapore: World
Scientific, 1994, 319 p.

\bibitem{Zi} A.A. Zinger. On independent samples from normal populations. (Russian).
Uspehi Matem. Nauk (N.S.), {\bf 6}, (1951), no. 5, 172--175.

\end{thebibliography}


\end{document}